\documentclass[a4paper,14pt]{article}
\usepackage[14pt]{extsizes}
\usepackage{mathtext}                    % Русский язык в формулах
\usepackage{cmap}                        % Поддержка поиска русских слов в PDF (pdflatex)
\usepackage[cp1251]{inputenc}            % Выбор языка и кодировки
\usepackage[english, russian]{babel}     % Настройки для русского языка
\usepackage[left=2.5cm,right=1cm,top=1cm,bottom=2cm]{geometry} % поля страницы
\usepackage{amssymb}
\usepackage{amsmath}

\newtheorem{lemma}{Лемма}
\newtheorem{theorem}{Теорема}

\newtheorem{corollary}{Следствие}

\newtheorem{definition}{Определение}

\newtheorem{proof}{Доказательство}
\newtheorem{remark}{Замечание}

\begin{document}

%\udk{517.983.246}

%\date{\\
%Исправленный вариант\\ }

%\author
%\address{Гомельский государственный университет им. Ф. Скорины,
%г.~Гомель}

\centerline{\textbf{Одно достаточное условие глобальной операторной монотонности}}

\centerline{\textbf{функций нескольких переменных}}

\centerline{А.\,Р.~Миротин}

\centerline{amirotin@yandex.ru}

Классическая теорема Лёвнера (см. \cite{Le1} --- \cite{BS}), описывающая   операторно-монотонные
функции одного переменного, обобщалась на случай нескольких  переменных многими авторами (см., например, работы \cite{AMY} --- \cite{H} и цитированную там литературу). В  частности, в  \cite{AMY} --- \cite{AMY3} дан критерий локальной операторной монотонности и указаны достаточные условия глобальной операторной монотонности функций нескольких  переменных.
Данное сообщение содержит два замечания относительно моего краткого сообщения  \cite{MZ}, также посвященного  достаточным условиям глобальной операторной монотонности. Первое из них обобщает и исправляет   основной результат указанного сообщения, во втором  выделен новый класс операторно-монотонных функций нескольких   переменных.

1. Следующее определение расширяет класс функций, рассмотренный в  \cite{MZ}.

\begin{definition}\label{definition1}
Скажем, что неположительная непрерывная на $(-\infty,0]^n$ функция
$\psi$ принадлежит классу $\mathcal{R}^-_n$, если
  существует такие числа $\lambda >0, \gamma\geq 0$ и
 положительная мера $\tau$ на $\Bbb{R}_+^n$, что при всех $z\in \mathbb{C}^n$ с ${\rm
Re}z_j>0\ (j=1,\ldots ,n)$ выполняется равенство
$$
\frac{1}{\lambda -\psi(-z)}=\gamma+S\tau(z),
$$
где
$$
S\tau(z)=\int\limits_{\Bbb{R}_+^n} \frac{d\tau(\xi)}
{(\xi_1+z_1)\ldots(\xi_n+z_n)}
$$
суть $n$-мерное преобразование Стилтьеса меры $\tau$.
\end{definition}

Далее нам понадобится также следующий класс функций.

\begin{definition}\label{definition2}
Пусть $U$ есть открытое подмножество $\mathbb{R}^n$, $\alpha\in \mathbb{R}$. Скажем, что непрерывная положительная функция  $f$, определенная на $U$, принадлежит классу $\mathcal{Q}^\alpha(U)$, если
  существуют такое число $\gamma$ и
 положительная мера $\tau$ на $\Bbb{R}_+^n$, что при всех $z\in U$  выполняется равенство
$$
\frac{1}{f(z)}=\gamma+S^\alpha\tau(z),\eqno(1)
$$
где
$$
S^\alpha\tau(z)=\int\limits_{\Bbb{R}_+^n} \frac{d\tau(\xi)}
{((\xi_1+z_1)\ldots(\xi_n+z_n))^\alpha}
$$
суть $n$-мерное обобщенное преобразование Стилтьеса  меры  $\tau$.
\end{definition}

Таким образом, если  $\psi\in \mathcal{R}^-_n$, то   $\lambda-\psi(-z)\in \mathcal{Q}^1((0,\infty)^n)$ при некотором  $\lambda>0$. Поэтому нижеследующая лемма 1 обобщает теорему 2 из  \cite{MZ}.
  Всюду ниже $\sigma(A)$ обозначает совместный спектр набора $A$ попарно коммутирующих самосопряженных операторов в комплексном гильбертовом пространстве (см., например, \cite{IZVR}), $E$ --- произведение спектральных мер операторов $A_j\ (j=1,\dots,n)$. Для любой непрерывной в окрестности  $\sigma(A)$ функции $f$, как известно, можно
определить оператор $f(A)$ формулой
$$
f(A)=\int\!\!\dots\!\!\int_{\sigma(A)}\!f(s)dE(s).
$$
При этом соответствие $f\mapsto f(A)$ однородно, аддитивно и мультипликативно (\cite[п. 111]{RN},  \cite[глава 6, \S 5]{BSol}).

\begin{lemma}\label{lemma1}
Пусть $f\in \mathcal{Q}^\alpha(U),\ U\subseteq (0,\infty)^n,\ 0\leq \alpha\leq 1$ и выполняется равенство $\mathrm{(1)}$. Для любого набора $A$ попарно коммутирующих ограниченных самосопряженных операторов в гильбертовом пространстве, таких, что  $\sigma(A)\subseteq U$, справедливо равенство

$$
f(A)^{-1}=
\gamma I+\int\limits_{\Bbb{R}_+^n}
\left(\prod_{j=1}^n(\xi_jI +A_j)\right)^{-\alpha} d\tau(\xi).
$$
\end{lemma}

\begin{proof} Интерес представляет лишь случай $\alpha>0$. Имеем
$$
 f(A)^{-1}=\int\limits_{\sigma(A)}\frac{1}{f(s)}dE(s)=\int\limits_{\sigma(A)}\left(\gamma+\int\limits_{\Bbb{R}_+^n} \frac{d\tau(\xi)}
{\left(\prod_{j=1}^n(\xi_j+s_j)\right)^\alpha}\right) dE(s)=
$$
$$
\gamma I+\int\limits_{\Bbb{R}_+^n} \int\limits_{\sigma(A)}\left(\prod_{j=1}^n(\xi_j+s_j)\right)^{-\alpha}dE(s)d\tau(\xi) =
\gamma I+\int\limits_{\Bbb{R}_+^n}
\left(\prod_{j=1}^n(\xi_jI +A_j)\right)^{-\alpha} d\tau(\xi).
$$
 Последнее равенство при $\alpha=1$ вытекает  из  аддитивности и мультипликативности соответствия $f\mapsto f(A)$.  При $0<\alpha<1$ его можно  обосновать с помощью формулы
$$
C^{-\alpha}=\frac{\sin(\alpha\pi)}{\pi}\int\limits_0^\infty t^{-\alpha}(tI+C)^{-1}dt\ (C>0,\ 0<\alpha<1).
$$
В самом деле, совместный спектр  $\sigma(A)$ совпадает с совместным аппроксимативным спектром набора $A$ \cite[теорема 2]{IZVR}, а потому обладает проекционным свойством \cite{SZ}. Следовательно, условие   $\sigma(A)\subseteq U$ влечет $\forall j \sigma(A_j)\subseteq (0,\infty)$. Таким образом, оператор $C:=\prod_{j=1}^n(\xi_jI +A_j)$  положительный вслед за $A_j$. Поэтому, дважды применяя упомянутую формулу для дробной степени, получаем  в силу теоремы Фубини
$$
\int\limits_{\sigma(A)}\left(\prod_{j=1}^n(\xi_j +s_j)\right)^{-\alpha}dE(s)=
\int\limits_{\sigma(A)}\frac{\sin(\alpha\pi)}{\pi}\int\limits_0^\infty t^{-\alpha}\left(t+\prod_{j=1}^n(\xi_j +s_j)\right)^{-1}dtdE(s)=
$$
$$
\frac{\sin(\alpha\pi)}{\pi}\int\limits_0^\infty t^{-\alpha}\int\limits_{\sigma(A)}\left(t+\prod_{j=1}^n(\xi_j +s_j)\right)^{-1}dE(s)dt=
$$
$$
\frac{\sin(\alpha\pi)}{\pi}\int\limits_0^\infty t^{-\alpha}\left(tI+\prod_{j=1}^n(\xi_jI +A_j)\right)^{-1}dt=
\left(\prod_{j=1}^n(\xi_jI +A_j)\right)^{-\alpha}.
$$
\end{proof}

 Теорема 1 сообщения \cite{MZ} должна выглядеть следующим образом:
\begin{theorem}\label{theorem1} Пусть функция $\psi$ принадлежит $\mathcal{R}^-_n$.
 Тогда
для  наборов $A=(A_1,\dots,A_n)$ и $A'=(A'_1,\dots,A'_n)$
ограниченных  неположительных операторов в
$H$ неравенства $A_j\leq A_j' \ (j=1,\dots, n)$ влекут неравенство
$\psi(A)\leq \psi(A')$, если операторы   $A_1,\dots,A'_n$  попарно коммутируют.
\end{theorem}

\begin{proof}
Воспользуемся  леммой 1
  и тем фактом, что  функция одного переменного  $z\mapsto -z^{-1}$  является операторно-монотонной (доказательство этого утверждения, независимое от
теоремы Лёвнера, см. в \cite{BB}, глава II, \S 47). В силу последнего утверждения  при всех $\xi_j\geq 0$   имеем $0<R(\xi_j ,A_j)\leq R(\xi_j, A_j')\ (j=1,\dots,n)$. Поскольку неравенства между положительными операторами, которые попарно коммутируют, можно почленно перемножать,
то и $0<\prod_{j=1}^nR(\xi_j ,A_j)\leq \prod_{j=1}^nR(\xi_j
,A_j')$. Тогда   $0<R(\lambda, \psi(A))\leq R(\lambda, \psi(A'))$
при некотором $\lambda>0$ по лемме 1, где положено $f(z)=\lambda-\psi(-z), \alpha=1$, а набор  $A$  заменен набором $-A$, и осталось
применить к этому неравенству функцию $-z^{-1}$.
\end{proof}

\begin{remark}
  Таким образом (при $n>1$), по сравнению с первоначальной формулировкой теоремы 1 мы добавляем требование коммутирования каждого из операторов $A_j$ с
каждым из операторов $A'_k\ (j, k=1,\dots,n)$. Без этого дополнительного требования теорема 1 работы \cite{MZ}  неверна. Действительно, легко видеть, что функция $\lambda-z_1z_2 (\lambda>0)$  принадлежит $\mathcal{R}^-_2$. В то же время, она не является глобально операторно-монотонной на $(-\infty,0)^2$ (см. определение 3 ниже), так как, например, не удовлетворяет необходимому условию глобальной операторной монотонности \cite[теорема 8.1]{AMY} (см. также \cite{AMY3}), а именно,  не принадлежит классу Лёвнера $\mathcal{L}(E)$, поскольку не отображает $\Pi^2$ в $\overline{\Pi}$, где $\Pi$ --- открытая верхняя полуплоскость, а $\overline{\Pi}$ --- ее замыкание.
\end{remark}

2.  Вместе с тем, можно модифицировать  рассуждения заметки \cite{MZ} и получить теорему, свободную от указанного дополнительного требования.

 \begin{definition}\label{definition3}
  \cite{AMY}. Пусть $U$ есть открытое подмножество $\mathbb{R}^n$, и  $f$  есть вещественнозначная непрерывная функция на $U$. Функция  $f$ называется (глобально) операторно-монотонной на $U$, если для любых наборов $A$ и $B$ из  $n$   попарно коммутирующих ограниченных самосопряженных операторов в гильбертовом пространстве, таких, что $\sigma(A)\cup\sigma(B)\subset U$,   из неравенств $A_j\leq B_j\ (j=1,\dots,n)$ следует неравенство $f(A) \leq f(B)$.
\end{definition}

Следующий результат для $n=2$ доказан в \cite{AMY} (см. также \cite{AMY3}).

\begin{lemma}\label{lemma2}
 Функция $(z_1\dots z_n)^\alpha$ является
 операторно-монотонной  на множестве $(0,\infty)^n$ при $ 0\leq \alpha\leq 1/n$.
\end{lemma}

 \begin{proof} Заметим, что композиция операторно-монотонных функций есть операторно-монотонная функция (при подходящем выборе областей определения). Так как по теореме Лёвнера (см., например, \cite[ с. 41]{RR}) функция одного переменного $z^\beta$ является  операторно-монотонной на $(0,\infty)$ при $0\leq\beta\leq 1$,  тождество  $(z_1\dots z_n)^\alpha=$ $((z_1\dots z_n)^{1/n})^{n\alpha}$ позволяет свести доказательство к случаю  $\alpha=1/n$. При $n=2$, как уже отмечалось, утверждение леммы доказано в \cite{AMY}, \cite{AMY3} (см. там пример 9.8). Если $n$ является степенью двойки, $n=2^l$, оно доказывается индукцией по  $l$ в силу тождества
$$
(z_1\dots z_{2^{l+1}})^{\frac{1}{2^{l+1}}}=\left((z_1\dots z_{2^l})^{\frac{1}{2^l}}(z_{2^l+1}\dots z_{2^{l+1}})^{\frac{1}{2^l}}\right)^{\frac{1}{2}}.
$$

Для доказательства нашего утверждения для любого $n$ применим теперь обратную индукцию. Если оно верно при некотором $n\geq 3$, то функция $(z_1\dots z_{n-1})^{1/n}$  будет
 операторно-монотонной  на $(0,\infty)^{n-1}$. Но тогда  операторно-монотонной  на $(0,\infty)^{n-1}$ будет и функция
 $$
 (z_1\dots z_{n-1}(z_1\dots z_{n-1})^{\frac{1}{n}})^{\frac{1}{n}}=(z_1\dots z_{n-1})^{\frac{1}{n}+\frac{1}{n^2}}.
 $$
Продолжая этот процесс, получим, что для любого  $k\in \mathbb{N}$ операторно-монотон\-ными  на $(0,\infty)^{n-1}$ будут  функции
$(z_1\dots z_{n-1})^{p_k}$, где  $p_k=1/n+\dots+1/n^k$, причем $p_k\to 1/(n-1)$ при $k\to\infty$.
Рассмотрим наборы $A$ и $B$ из  $n-1$   попарно коммутирующих ограниченных самосопряженных операторов, такие, что $\sigma(A)\cup\sigma(B)\subset (0,\infty)^{n-1}$   и $A_j\leq B_j\ (j=1,\dots,n-1)$. Тогда $(A_1\dots A_{n-1})^{p_k}\leq (B_1\dots B_{n-1})^{p_k}$, что в силу известного свойства дробных степеней положительных операторов в пределе дает требуемое неравенство.
\end{proof}

\begin{corollary}
 Функция  $\sum_{j=1}^nc_j(z_1\dots z_{j})^{\alpha_j}$ при $ 0\leq \alpha_j\leq 1/n,\ c_j\geq 0$ является операторно-монотонной  на
$(0,\infty)^n$.
\end{corollary}

Следующая теорема обобщает лемму 2.

\begin{theorem}\label{theorem2}
Если $f\in \mathcal{Q}^\alpha((0,\infty)^n)$ при некотором  $\alpha$, $0\leq \alpha\leq 1/n$, то $f$ есть  операторно-монотонная функция на $(0,\infty)^n$.
\end{theorem}

 \begin{proof} Пусть  $A$ и $B$ есть  наборы из  $n$   попарно коммутирующих ограниченных самосопряженных операторов, причем $\sigma(A)\cup\sigma(B)\subset (0,\infty)^{n}$   и $A_j\leq B_j\ (j=1,\dots,n)$. Если выполняется равенство (1), то по лемме 1
$$
f(A)^{-1}=
\gamma I+\int\limits_{\Bbb{R}_+^n}
\left(\left(\prod_{j=1}^n(\xi_jI +A_j)\right)^{\alpha}\right)^{-1}d\tau(\xi).
$$
В силу леммы 2
$$
\left(\prod_{j=1}^n(\xi_jI +A_j)\right)^{\alpha}\leq \left(\prod_{j=1}^n(\xi_jI +B_j)\right)^{\alpha}.
 $$
 Поэтому операторная монотонность функции $-z^{-1}$ на $(0,\infty)$  влечет неравенство $f(B)^{-1}\leq f(A)^{-1}$, а вслед за ним и неравенство $f(A)\leq f(B)$, что и требовалось доказать.
\end{proof}

\begin{corollary}
 Функции вида  $f(z)=1/(\gamma+\prod_{j=1}^nf_j(z_j))$, где
 $$
f_j(x)= \int\limits_{0}^\infty \frac{d\tau_j(s)}
{(s+x)^\alpha},\quad  0\leq \alpha\leq 1/n
 $$
 есть обобщенное преобразование Стилтьеса положительной меры $\tau_j$,
 являются оператор\-но-монотонными  на $(0,\infty)^n$.
\end{corollary}

\begin{remark}
    Лемма 2, а вместе с ней  и теорема 2, неверны при $\alpha\notin [0,1/n]$, так как при таких $\alpha$  функция $(z_1\dots z_n)^\alpha$  не удовлетворяет уже упоминавшемуся в замечании 1 необходимому условию глобальной операторной монотонности из \cite[теорема 8.1]{AMY} (см. также \cite{AMY3}), а именно,  не отображает $\Pi^n$ в $\overline{\Pi}$ (как и выше, $\Pi$ --- открытая верхняя полуплоскость,  $\overline{\Pi}$ --- ее замыкание).
\end{remark}

\end{document}